# Non-Negative Integer-Valued Semi-Selfsimilar Processes

S Satheesh

NEELOLPALAM, S. N. Park Road
Trichur – 680 004, **India.**

*ssatheesh1963@yahoo.co.in*

and

E Sandhya

Department of Statistics, Prajyoti Niketan College
Pudukkad, Trichur – 680 301, **India.**

*esandhya@hotmail.com*

***Abstract*:** Non-negative integer-valued semi-selfsimilar processes are introduced. Levy processes in this class are characterized. Its relation to an AR(1) scheme is derived.

1. **Introduction**

Semi-selfsimilar (SSS) processes are stochastic processes that are invariant in distribution under suitable scaling of time and space. That is; a process $\{X(t), t \geq 0\}$ is SSS if for some $a>0$ there exists an $H>0$ such that,

$$\{X(at)\} \stackrel{d}{=} \{a^H X(t)\}. \qquad (1.1)$$

If the above relation holds for any $a>0$, then $\{X(t)\}$ is selfsimilar (SS). Here $a$ is called the epoch and $H$ the exponent of the SSS process and we write $\{X(t)\}$ is $(a,H)$-SSS. Clearly for the SS process we need specify $H$ only and hence we write $\{X(t)\}$ is $H$-SS. $\{X(t), t \geq 0\}$ is a Lévy process if it has stationary and independent increments and $X(\text{o}) = 0$. Recall that $\{X(t)\}$ is Lévy *iff* $X(1)$ is infinitely divisible (ID). Maejima and Sato (1999) showed that a Lévy process $\{X(t)\}$ is SSS (SS) *iff* the distribution of $X(1)$ is semi-stable (stable). For more on this see also Embrechts and Maejima (2000).

The notion of semi-stable laws is by Lévy. A distribution on **R** is semi-stable($a,b$) if its characteristic functions (CF) $f(u)$ satisfies $f(u) \neq 0$ and $\{f(u)\}^a = f(bu)$, $\forall u \in \mathbf{R}$ and some $b \in (0,1) \cup (1,\infty)$. Then there exists a unique $\alpha \in (0,2]$ satisfying $b = a^{1/\alpha}$. Similarly, a distribution on $\mathbf{R}^+ = [0,\infty)$ is semi-stable if its Laplace transforms (LT) $\varphi(s)$ satisfies $\{\varphi(s)\}^a = \varphi(bs)$, $\forall s > 0$ and some $b > 0$. In this case we have a unique $\alpha \in (0,1]$ satisfying $b = a^{1/\alpha}$. Semi-stable laws on $I_0$ were described in terms of PGFs by Satheesh and Nair (2002) and Bouzar (2004).

Of late there has been an increasing interest in describing and discussing the integer-valued analogues of classical distributions *viz.*, stable, semi-stable, semi-selfdecomposable (semi-SD), geometrically ID, Laplace *etc*. See, eg. Aly and Bouzar (2000), Bouzar (2004, 2006), Bouzar and Jayakumar (2006), Inusah and Kozubowski (2006), Kozubowski and Inusah (2006) and Satheesh and Nair (2002).

The aim of this note is to extend the notion of semi-selfsimilarity to integer-valued processes. Clearly, equation (1.1) does not (cannot) include the integer-valued case. Instead, we rely on the idea of binomial thinning to arrive at the right definitions. In the next section we describe non-negative integer-valued ($I_0$-valued) SSS processes and characterize $I_0$-valued SSS Lévy processes in terms of an $I_0$-valued first order autoregressive (AR(1)) series. This characterization is the $I_0$-valued analogue of the discussion in Satheesh and Sandhya (2007).

## 2. Results

The description in (1.1) cannot hold good for an integer-valued process, as the state space of the process on the RHS will be different. Thus it is clear that the scaling of the state space as given in equation (1.1) does not hold for integer-valued *r.v*s. Hence we need a different formulation for integer-valued SSS processes. This formulation is based on the idea of binomial thinning due to Steutel and van Harn (1979). The binomial thinning operation "$\otimes$" on an $I_0$-valued r.v $X$ is defined as: $b \otimes X = \sum_{i=1}^{X} Z_i$ for some i.i.d Bernoulli r.vs $\{Z_i\}$ independent of $X$ with $P\{Z_i=1\}= b$. They viewed this operation as the $I_0$-valued analogue of scalar multiplication of continuous r.vs. A justification for using $b \otimes X$ as the $I_0$-valued analogue of $bX$ in the continuous set up was given by Satheesh and Nair (2002). It may be noted that if $P(s)$ is the probability generating function (PGF) of $X$ then that of $b \otimes X$ is $P(1-b+bs)$.

**Definition.2.1** An $I_0$-valued r.v $X$ is discrete semi-stable$(a,b)$ *if* its PGF $P(s)$ satisfies $P(s) \neq 0$ and $\{P(s)\}^a = P(1-b+bs)$, $\forall\ 0 \leq s \leq 1$, some $0<b<1$ and $b = a^{1/\alpha}$, $\alpha \in (0,1]$.

We now give an example each of an $\mathbf{R}^+$-valued and $I_0$-valued semi-stable law.

**Examples.** By Schoenberg (1938) there is a one-to-one correspondence between a real CF and a LT. Pillai (1985) showed that

$$\omega(u) = \frac{1}{1+|u|^\alpha\ \{1 - A\cos[k\log(|u|)]\}},\ u \in \mathbf{R},\ k = \frac{-2\pi}{\log(b)},\ 0<b<1,\ 0<A<1,\ 0<\alpha \leq 2$$

is a Polya type CF. From this CF we have $\psi(u) = |u|^\alpha\ \{1 - A\cos[k\log(|u|)]\}$ as an example of the function satisfying $a\psi(u) = \psi(bu)$, $u \in \mathbf{R}$, $b^\alpha = a$, $\alpha \in (0,2]$, and consequently a description of semi-stable CFs of the form $exp\{-\psi(u)\}$. Noticing that this is a real CF we have the LT of an $\mathbf{R}^+$-valued semi-stable law as;

$$\varphi(s) = \exp\{-s^\alpha\{1 - A\cos[k\log(s)]\}\},\ s>0,\ k = \frac{-2\pi}{\log(b)},\ 0<b<1,\ 0<A<1,\ \alpha \in (0,1].$$

Now if $\varphi(s)$ is a LT, then $P(s) = \varphi(1-s)$, $0 \leq s \leq 1$ is a PGF, being that of a mixture of Poisson laws, the mixing distribution having LT $\varphi$. Hence for $0 \leq s \leq 1$;

$$\exp\{-(1-s)^{\alpha}\{1-A\cos[k\log(1-s)]\}\}, \ k = \frac{-2\pi}{\log(b)}, \ 0<b<1, \ 0<A<1, \ \alpha \in (0,1],$$

is the PGF of an $I_0$-valued semi-stable law. For another example see Bouzar (2004).

**Definition.2.2** An $I_0$-valued r.v $X$ with PGF $P(s)$ is discrete semi-SD($b$) *if* there exists a PGF $P_o(s)$ that is ID such that $P(s) = P(1-b+bs) P_o(s), \forall \ 0 \leq s \leq 1$ and some $0<b<1$.

Now the $I_0$-valued analogue of SSS process is as follows.

**Definition.2.3** An $I_0$-valued process $\{X(t), t \geq 0\}$ is SSS *if* for some $0<a<1$ there exists an $H>0$ such that

$$\{X(at)\} \stackrel{d}{=} \{a^H \otimes X(t)\}.$$

If the above relation holds for any $0<a<1$, then $\{X(t)\}$ is SS. Also, as in the continuous case we call $a$ the epoch and $H$ the exponent of the SSS (or SS) process and we write $\{X(t)\}$ is $(a,H)$-SSS, and the SS process as $\{X(t)\}$ is $H$-SS.

**Note.** Here the notions of semi-stability and semi-selfsimilarity are considered in the strict sense only.

Notice that if $X(t)$ is an $I_0$-valued Lévy process and $P(s)$ is the PGF of $X(1)$, then that of $X(t)$ is $P(s)^t$. Next we give a corollary to theorem.4.1 in Maejima and Sato (1999) and follow it with its $I_0$-valued analogue.

**Theorem.2.1** An $\mathbf{R}^+$-valued Lévy process $\{X(t), t \geq 0\}$ is $(a, \frac{1}{\alpha})$-SSS ($\frac{1}{\alpha}$-SS) *iff* the distribution of $X(1)$ is $\mathbf{R}^+$-valued semi-stable($a,b$) ($\alpha$-stable), $\alpha \in (0,1]$.

**Theorem.2.2** An $I_0$-valued Lévy process $\{X(t), t \geq 0\}$ is $(a, \frac{1}{\alpha})$-SSS *iff* the distribution of $X(1)$ is $I_0$-valued semi-stable($a,b$), $\alpha \in (0,1]$.

*Proof.* If $X(t)$ is $I_0$-valued Lévy and $X(1)$ is $I_0$-valued semi-stable($a,b$) then for each $t \geq 0$,

$$\{P(s)\}^{at} = \{P(1-b+bs)\}^{t}, \forall\ 0 \leq s \leq 1,\ \text{some}\ 0<b<1\ \text{and}\ b = a^{1/\alpha}.\ \text{That is};$$

$$X(at) \stackrel{d}{=} a^{1/\alpha} \otimes X(t)\ \text{for each}\ t \geq 0,\ \text{since}\ b = a^{1/\alpha}.$$

Thus $\{X(t)\}$ is $I_0$-valued $(a, \frac{1}{\alpha})$-SSS.

Conversely, if a Lévy process $\{X(t)\}$ is $I_0$-valued $(a, \frac{1}{\alpha})$-SSS then

$$X(at) \stackrel{d}{=} a^{1/\alpha} \otimes X(t),\ \text{and hence in terms of PGFs at}\ t=1,$$

$$\{P(s)\}^{a} = \{P(1 - a^{1/\alpha} + a^{1/\alpha}s)\}.$$

Hence the PGF $P(s)$ is semi-stable$(a, a^{1/\alpha})$.

**Remark.2.1** The CF $f(u)$ of a semi-stable$(a,b)$ law can also be equivalently written as $\{f(\frac{u}{b})\}^{a} = f(u)$. Since $b = a^{1/\alpha}$, the Lévy process version of this is: $\frac{1}{b}X(t) \stackrel{d}{=} X(\frac{t}{a})$ or $X(\frac{t}{a}) \stackrel{d}{=} (\frac{1}{a})^{1/\alpha} X(t)$. That is, a semi-stable$(a,b)$ law gives rise to a $(\frac{1}{a}, \frac{1}{\alpha})$-SSS Lévy process also. Thus in the real or $\mathbf{R}^{+}$-valued cases, SSS Lévy process corresponding to the same semi-stable$(a,b)$ law has two epochs, $a$ and $\frac{1}{a}$. However, in the $I_0$-valued setup the epoch can only be $a$, since we need $a<1$, $a^{1/\alpha}$ being a Bernoulli probability.

**Theorem.2.3** An $I_0$-valued semi-stable$(a,b)$ law is semi-SD$(b)$.

*Proof.* Let $P(s)$ be the PGF of the semi-stable$(a,b)$ law. Then

$$P(s) = \{P(1-b+bs)\}^{1/a},\ \forall\ 0 \leq s \leq 1,\ \text{some}\ 0<b<1\ \text{and}\ b = a^{1/\alpha}.$$

$$= P(1-b+bs)\{P(1-b+bs)\}^{\frac{1}{a}-1}.$$

Here the second factor is also ID since $P(s)$ is ID (Bouzar (2007)) and hence by definition.2.2 $P(s)$ is semi-SD$(b)$. For another proof, see Bouzar (2007).

Thus by Satheesh and Sandhya (2005) $I_0$-valued semi-stable($a,b$) family of laws qualify to model stationary additive $I_0$-valued AR(1) schemes described by the sequence of r.v $\{X_n\}$, if there exists an innovation sequence $\{\varepsilon_n\}$ of *i.i.d* r.vs satisfying;

$$X_n = b \otimes X_{n-1} + \varepsilon_n, \; n>0 \text{ integer and some } 0<b<1. \tag{2.1}$$

Rather than just prescribing $\varepsilon_n$ to have the PGF $\{P(1-b+bs)\}^{\frac{1}{a}-1}$, we take a different look at this. Here we make use of the SSS processes we briefly discussed.

**Theorem.2.4** Let $\{Z(t), t \geq 0\}$ be an $I_0$-valued Levy process and in (2.1) let $X_0 \stackrel{d}{=} Z(1)$ and $\varepsilon_n \stackrel{d}{=} b \otimes Z(b^{-\alpha}-1)$, $\forall n$, $0<b<1$, $0<\alpha\leq 1$. Then (2.1) is stationary with $I_0$-valued semi-stable($b^\alpha,b$) marginals *if* $\{Z(t)\}$ is $(b^\alpha, \frac{1}{\alpha})$-SSS. Conversely, the marginals of (2.1) are $I_0$-valued semi-stable($b^\alpha,b$) and $\{Z(t)\}$ $(b^\alpha, \frac{1}{\alpha})$-SSS *if* (2.1) is stationary.

**Proof.** Let the PGF of $Z(1)$ is $P(s)$, then that of $b \otimes Z(b^{-\alpha}-1)$ is $\{P(1-b+bs)\}^{b^{-\alpha}-1}$. If $\{Z(t)\}$ is $(b^\alpha, \frac{1}{\alpha})$-SSS then,

$$Z(b^\alpha t) \stackrel{d}{=} b \otimes Z(t) \text{ or } \{P(s)\}^{b^\alpha} = P(1-b+bs), \text{ and } Z(1) \text{ is semi-stable}(b^\alpha,b).$$

Now under the given assumptions at $n=1$ in (2.1), the PGF of $X_1$ is

$$P_1(s) = P(1-b+bs)\{P(1-b+bs)\}^{b^{-\alpha}-1} = \{P(1-b+bs)\}^{b^{-\alpha}} = P(s),$$

which is same as that of $X_0$. Hence on iteration (2.1) is stationary with $I_0$-valued semi-stable($b^\alpha,b$) marginals.

Conversely, let (2.1) is stationary and $P(s)$ is the PGF of $X_0$. Then at $n=1$,

$$P(s) = P(1-b+bs)\{P(1-b+bs)\}^{b^{-\alpha}-1} = \{P(1-b+bs)\}^{b^{-\alpha}}.$$

Hence the marginals and $Z(1)$ are semi-stable($b^\alpha,b$) and consequently $\{Z(t)\}$ is $(b^\alpha, \frac{1}{\alpha})$-SSS.

**Concluding Remarks.** As corollaries to theorems 2.2 and 2.4 we have (in the $I_0$-valued setup): If $\{X(t)\}$ is stable Levy then $\{X(t)\}$ is $\frac{1}{\alpha}$-SS. Similarly, (2.1) is stationary with stable marginals if $\{Z(t)\}$ is $\frac{1}{\alpha}$-SS. Conversely, $\{Z(t)\}$ is SS and the marginals are stable if (2.1) is stationary. This is because if the condition describing semi-stability and semi-selfsimilarity is true for two reals $a_1$ and $a_2$ such that $\log(a_1)/\log(a_2)$ is irrational then it describes stability and selfsimilarity, see Embrechts and Maejima (2000).

**References.**